\documentclass{amsart}
\usepackage{amsmath}
\usepackage{graphicx}
\usepackage{amsfonts}
\usepackage{amssymb}
\newtheorem{theorem}{Theorem}
\theoremstyle{plain}

\newtheorem{corollary}{Corollary}

\newtheorem{definition}{Definition}
\newtheorem{example}{Example}

\newtheorem{lemma}{Lemma}

\newtheorem{proposition}{Proposition}
\newtheorem{remark}{Remark}

\numberwithin{equation}{section}

\begin{document}
\title[Characterization of optimal Transport Plans]{Characterization of optimal Transport Plans for the Monge-Kantorovich-Problem}
\author{Walter Schachermayer and Josef Teichmann}
\address{Technical University Vienna, Wiedner Hauptstrasse 8--10, A-1040 Vienna, Austria}
\thanks{Financial support from the Austrian Science Fund (FWF) under
  grant P 15889, from the Vienna Science Foundation (WWTF) under grant
  MA13, from the European Union under grant HPRN-CT-2002-00281 is
  gratefully acknowledged. Furthermore this work was financially
  supported by the Christian Doppler Research Association (CDG). The
  authors gratefully acknowledge a fruitful collaboration and
  continued support by Bank Austria through CDG} 
\begin{abstract}
We prove that $c$-cyclically monotone transport plans $\pi$ optimize
the Monge-Kantorovich transportation problem under an additional
measurability condition. This measurability condition is always
satisfied for finitely valued, lower semi-continuous cost functions. In
particular, this yields a positive answer to Problem 2.25 in C.
Villani's book. We emphasize that we do not need any regularity
conditions as were imposed in the previous literature.
\end{abstract}\maketitle

We consider the \emph{Monge-Kantorovich optimal transport problem }$(\mu
,\nu,c)$ for Borel measures $\mu,\nu$ on polish spaces $X,Y$ and a lower
semi-continuous cost function $c:X\times Y\to\mathbb{R}_{\geq0}
\cup\{+ \infty\}$ (see C.~Villani's beautiful book \cite{vil:03} for
all necessary details): we define the set of \emph{transport plans}
$\Pi(\mu,\nu)$ as the set of probability measures $\pi$ with marginal
$\mu$ on $X$, and marginal $\nu$ on $Y$, respectively. Furthermore, we
define $\Phi(\mu,\nu)$ as the set of pairs $(\phi,\psi)$ of Borel
functions $\phi:X\to \mathbb{R\cup\{-\infty\}}$ and
$\psi:Y\to\mathbb{R\cup\{-\infty\}}$ with
$\phi\in\mathcal{L}^{1}(\mu)$ and $\psi\in\mathcal{L}^{1}(\nu)$, such
that
\[
\phi(x)+\psi(y)\leq c(x,y)
\]
for all $(x,y)\in X\times Y$. The Monge-Kantorovich problem is to minimize the
cost
\[
I(\pi)=\int_{X\times Y}c(x,y)\pi(dx,dy)
\]
for transport plans $\pi\in\Pi(\mu,\nu)$. Dually to the Monge-Kantorovich
problem, we maximize
\[
J(\phi,\psi)=\int_{X}\phi(x)\mu(dx)+\int_{Y}\psi(y)\nu(dy)
\]
over $(\phi,\psi)\in\Phi(\mu,\nu)$.

It is well-known that under the previous conditions, the equality
\[
\inf_{\Pi(\mu,\nu)}I(\pi)=\sup_{\Phi(\mu,\nu)}J(\phi,\psi)
\]
holds true, and the infimum is in fact a minimum, if $\inf_{\Pi(\mu,\nu)}
I(\pi)<\infty$ (see \cite{vil:03}, Th.1.3). We know furthermore that
optimizers $\widehat{\pi}$ with $I(\widehat{\pi})<\infty$ are concentrated on
$c$-monotone (see the following definition) Borel sets $\Gamma\subset X\times
Y$. These assertions are proved for instance in \cite[ch.2]{vil:03} or in
\cite[Th. 3.2]{ambpra:02}.

\begin{definition}
A set $\Gamma\subset X\times Y$ is called $c$\emph{-monotone} if for
all $n\geq1$, for all points $(x_{i},y_{i})\in\Gamma$, $i=1,\dots,n$,
and all
permutations $\sigma\in\mathfrak{S}_{n}$
\[
\sum_{i=1}^{n}c(x_{i},y_{i})\leq\sum_{i=1}^{n}c(x_{i},y_{\sigma(i)})
\]
holds true. A transport plan $\pi\in\Pi(\mu,\nu)$ is called $c$
\emph{-monotone} if there is a $c$-monotone Borel set $\Gamma\subset X\times
Y$, where $\pi$ is concentrated on, i.e. $\pi(\Gamma)=1$.
\end{definition}

In \cite{vil:03} it was posed as an open problem whether $c$-monotone
transport plans $\pi\in\Pi(\mu,\nu)$ are necessarily optimizers of the
Monge-Kantorovich problem (see Problem 2.25 in \cite{vil:03} and the
references therein). Due to \cite{ambpra:02} it is known that this
conjecture holds true under a certain moment conditions on the
probability measures $\mu,\nu$ with respect to the cost function $c$.

Crucial for these considerations is R\"{u}schendorf's Theorem (see
\cite{rue:95} and \cite[ch.2]{vil:03}), which relates $c$-concave functions
and $c$-monotone sets $\Gamma$. Given a function $\phi:X\to
\mathbb{R}\cup\{-\infty\}$ and a lower semi-continuous cost function
$c:X\times Y\to\mathbb{R}_{\geq0}\cup\{+\infty\}$. Then $\psi$ is
called $c$-concave on $Y$, if there is a function $\phi:X\to
\mathbb{R}\cup\{-\infty\}$, such that
\[
\psi(y)=\inf_{x\in X}(c(x,y)-\phi(x))
\]
for $y\in Y$. Accordingly we define $c$-concave functions $\phi$ on $X$.

Let $\phi:X\to\mathbb{R}\cup\{-\infty\}$ be a function not identical
to $- \infty$, then
\begin{equation}
\phi^{c}(y)=\inf_{x\in X}(c(x,y)-\phi(x)) \tag{1}\label{c conjugate}
\end{equation}
is called the $c$-transform of $\phi$, which is obviously a $c$-concave
function on $Y$. We have that $\phi$ is $c$-concave on $X$ if and only if
$\phi=\phi^{cc}$. The $c$-superdifferential $\partial_{c}\phi\subset X\times
Y$ of a $c$-concave function $\phi$ is defined as the set of all pairs
$(x,y)\in X\times Y$ such that for all $z\in X$,%
\[
\phi(z)\leq\phi(x)+(c(z,y)-c(x,y)).
\]
R\"{u}schendorf's Theorem finally states that a set $\Gamma\subset X\times Y$
is $c$-monotone if and only if there is $c$-concave function $\phi$ such that
$\Gamma\subset\partial_{c}\phi$ (see \cite[ch.2]{vil:03} and \cite{rue:95} for
further details and references).

We cite the Theorem from \cite{ambpra:02}, which inspired our subsequent considerations.

\begin{theorem}
Let $X,Y$ be polish spaces and $\mu,\nu$ be Borel probability measures
thereon. Assume that $c$ is a finitely-valued, lower semi-continuous cost
function and that $\pi$ is $c$-monotone. Furthermore we assume that
\begin{align*}
\nu\left(  \{y,\int_{X}c(x,y)\mu(dx)<\infty\}\right)   &  >0,\\
\mu\left(  \{x,\int_{Y}c(x,y)\nu(dy)<\infty\}\right)   &  >0.
\end{align*}
Then $\pi$ is an optimizer and there exists a dual optimizer $(\phi,\psi
)\in\mathcal{L}^{1}(\mu)\times\mathcal{L}^{1}(\nu)$, where $\phi$ is a
$c$-concave Borel function, $\psi$ is a Borel function, $\psi=\phi^{c}$ almost
surely with respect to $\nu$ and
\[
\phi(x)+\psi(y)\leq c(x,y)
\]
for $(x,y)\in X\times Y$, where equality holds $\pi$-almost surely.
\end{theorem}

We shall prove the same assertion without assuming that $c$ is finitely
valued and without any assumptions on the existence of moments of $c$.
Here we do assume a measurability condition related to $c$-monotone
transport plans, namely the existence of a Borel-measurable
``hull'' for a sequence of dual optimizers. We shall prove (Proposition \ref{existence} below) that this measurability condition always holds true for \emph{finitely valued}, lower semi-continuous cost functions $c$ thus covering in
particular the setting of problem 2.25 in \cite{vil:03}.

\begin{definition}
Given polish spaces $X,Y$, a lower semi-continuous cost function $c:X\times
Y\to\mathbb{R}_{\geq0}\cup\{+\infty\}$, we call a Borel set
$\Gamma\subset X\times Y$ \emph{strongly }$c$\emph{-monotone}, if there exist
Borel functions $\phi:X\to\mathbb{R}\cup\{-\infty\}$ and
$\psi:Y\to\mathbb{R}\cup\{-\infty\}$ such that%
\[
\phi(x)+\psi(y)\leq c(x,y),
\]
for all $(x,y)\in X\times Y$ and equality holds if $(x,y)\in\Gamma$. We call a
transport plan $\pi$ \emph{strongly }$c$\emph{-monotone} if there exists a
strongly $c$-monotone Borel set $\Gamma\subset X\times Y$, on which $\pi$ is
concentrated, i.e. $\pi(\Gamma)=1$. \label{strongly c monotone}
\end{definition}

\begin{remark}
A seemingly innocent, but important point is the Borel measurability of $\phi$
and $\psi$.
\end{remark}

\begin{remark}
Notice that we do neither assume the Borel-measurable functions $\phi$,
$\psi$ to be $c$-concave nor to be conjugate.
\end{remark}

\begin{remark}
A strongly $c$-monotone set $\Gamma\subset X\times Y$ is $c$-monotone, since
\begin{align*}
\sum_{i=1}^{n}c(x_{i},y_{i})  &  =\sum_{i=1}^{n}(\phi(x_{i})+\psi(y_{i}))\\
&  =\sum_{i=1}^{n}(\phi(x_{i})+\psi(y_{\sigma(i)}))\\
&  \leq\sum_{i=1}^{n}c(x_{i},y_{\sigma(i)}),
\end{align*}
for all $(x_{i},y_{i})_{i=1,\dots,n}\subset\Gamma$ and permutations $\sigma
\in\mathfrak{S}_{n}$ for $n\geq1$.\label{strong implies usual}
\end{remark}

The following Proposition, which is based on the proof of Theorem 3.2 in
\cite{ambpra:02}, shows that \textbf{all }$c$\textbf{-monotone transport plans
are strongly }$c$\textbf{-monotone} provided that $c$ is finitely valued.

\begin{proposition}
Given polish spaces $X,Y$, a finitely valued, lower semi-continuous cost
function $c:X\times Y\to\mathbb{R}_{\geq0}$ and a $c$-monotone
transport plan $\pi$. Then there exist Borel-measurable functions $\phi,\psi$
on $X,Y$, respectively, taking values in $\mathbb{R}\cup\{-\infty\}$, such
that $\phi$ is $c$-concave, the $c$-superdifferential $\partial_{c}\phi$
contains a set where $\pi$ is concentrated, $\psi$ coincides $\nu$-almost
surely with $\phi^{c}$, the inequality
\[
\phi(x)+\psi(y)\leq c(x,y)
\]
holds true for all $(x,y)\in X\times Y$, and equality holds true $\pi$-almost
surely. In particular $\pi$ is strongly $c$-monotone.\label{existence}
\end{proposition}

\begin{proof}
In order to apply a construction of \cite{ambpra:02} we have to show that
there is a $c$-monotone Borel set $\Gamma=\cup_{k\geq1}\Gamma_{k}$ with
$\Gamma_{k}$ compact, $c|_{\Gamma_{k}}$ continuous for $k\geq1$ and
$\pi(\Gamma)=1$. Take any $c$-monotone Borel set $\Gamma^{\prime}$ with
$\pi(\Gamma^{\prime})=1$, which exists by assumption. By Egorov's Theorem and
the inner regularity on polish spaces of the transport plan $\pi$ (see
\cite{kal:01} for both results), we know that there are compact subsets
$L_{k}\subset X\times Y$, such that $\pi(X\times Y\setminus L_{k})<\frac{1}
{k}$ for $k\geq1$ and $c|_{L_{k}}$ is continuous for $k\geq1$. Consequently we
can define compact subsets $\Gamma_{k}:=\overline{\Gamma^{\prime}\cap L_{k}
}\subset L_{k}\subset X\times Y$ for $k\geq1$. In particular $c|_{\Gamma_{k}}$
is continuous for $k\geq1$ and the Borel set $\Gamma:=\cup_{k\geq1}\Gamma_{k}$
satisfies $\Gamma\supset\Gamma^{\prime}$ by construction. Since $c|_{\Gamma
_{k}}$ is continuous and $\Gamma^{\prime}\cap L_{k}$ is $c$-monotone by
assumption, the $c$-monotonicity extends to the closure $\Gamma_{k}
:=\overline{\Gamma^{\prime}\cap L_{k}}$ for $k\geq1$. A union of $c$-monotone
sets is $c$-monotone, so $\Gamma$ is $c$-monotone. Whence we have the claim
and can apply \cite{ambpra:02}, Th. 3.2 (step 1 of the proof) in order to
construct -- in the spirit of R\"{u}schendorf's Theorem -- a Borel-measurable
function $\phi:X\to\mathbb{R}\cup\{-\infty\}$ with $\Gamma
\subset\partial_{c}\phi$.

We then define $\phi^{c}$ via (\ref{c conjugate}) and obtain by applying again
\cite{ambpra:02}, Th. 3.2 (step 2 of the proof), that $\phi^{c}$ is $\nu
$-measurable. Setting the values of $\phi^{c}$ to $-\infty$ on a $\nu
$-negligible set we can define a Borel measurable function $\psi$ on $Y$ such
that%
\[
\phi(x)+\psi(y)\leq c(x,y)
\]
holds true for all $(x,y)\in X\times Y$. If $(x,y)\in\Gamma\subset\partial
_{c}\phi$, we know from the definition of the superdifferential of $\phi$ that
\[
\phi(z)\leq\phi(x)+(c(z,y)-c(x,y))
\]
for all $z\in X$, hence
\[
c(x,y)-\phi(x)\leq c(z,y)-\phi(z)
\]
for all $z$, so
\[
c(x,y)-\phi(x)\leq\psi(y),
\]
for $(x,y)\in\Gamma$. Hence equality in
\[
\phi(x)+\psi(y)\leq c(x,y)
\]
holds true $\pi$-almost surely, since $\pi(\Gamma)=1$.
\end{proof}

Next we show that \textbf{optimal transport plans are strongly }
$c$\textbf{-monotone} provided that $c$ is $\mu\otimes\nu$-almost surely
finitely valued. This result suggests to work with strongly $c$-monotone
transport plans instead of $c$-monotone ones. For the proof we need the
following Lemma, which generalizes results in \cite{rue:95} and \cite{sve:87}
(see in particular the references therein).

\begin{lemma}
Let $(\Omega_{i},\mathcal{F}_{i},\mu_{i})$ be probability spaces and let $\phi
_{n}:\Omega_{1}\to\mathbb{R}$ and $\psi_{n}:\Omega_{2}\to
\mathbb{R}$ be measurable functions for $n\geq1$. Assume that
\[
\xi_{n}(\omega_{1},\omega_{2})=\phi_{n}(\omega_{1})+\psi_{n}(\omega_{2})
\]
converges in $\mu_1 \otimes \mu_2 $-probability to $\xi:\Omega_{1}\times
\Omega_{2}\to\mathbb{R\cup\{-\infty\}}$. Then there exist real numbers
$(r_{n})_{n\geq1}$ and measurable functions $\phi:\Omega_{1}\to
\mathbb{R\cup\{-\infty\}}$, $\psi:\Omega_{2}\to\mathbb{R\cup
\{-\infty\}}$, such that -- along a fixed subsequence  -- ${(\phi_{n}+r_{n})}_{n \geq 0} $ converges  in probability to $ \phi $ and ${(\psi_{n} -r_{n})}_{n \geq 0} $ to $ \psi $. Furthermore,
\[
\xi(\omega_{1},\omega_{2})=\phi(\omega_{1})+\psi(\omega_{2})
\]
$\mu_{1}\otimes\mu_{2}$-almost surely.\label{convergence}
\end{lemma}

\begin{proof}
We choose a complete, bounded metric $d$ on $\mathbb{R\cup\{-\infty\}}$. Then
$\xi_{n}\to\xi$ in probability as $n\to\infty$ is equivalent
to $E(d(\xi_{n},\xi))\to0$ as $n\to\infty$. By Fubini's
theorem the real-valued function
\[
\omega_{2}\mapsto\int_{\Omega_{1}}d(\xi_{n}(\omega_{1},\omega_{2}),\xi
(\omega_{1},\omega_{2}))\mu_{1}(d\omega_{1})
\]
is $\mu_{2}$-almost surely well-defined and converges in $\mu_{2}$-probability
to $0$ as $n\to\infty$ by assumption.
 
We shall now distinguish two cases, either $\xi$ is finitely valued
with positive probability, or $\mu_1 \otimes
\mu_2$-almost surely $\xi=-\infty$. First we assume that there is $A\subset\Omega_{1}\times\Omega_{2}$ such that $\xi|_{A}> -\infty$ and $(\mu_{1}\otimes\mu
_{2})(A)>0$. Then there is $\widetilde{\omega_{2}}\in\Omega_{2}$
such that
\[
\xi_{n}(\omega_{1},\widetilde{\omega_{2}})=\phi_{n}(\omega_{1})+\psi
_{n}(\widetilde{\omega_{2}})\to\phi(\omega_{1}):=\xi(\omega
_{1},\widetilde{\omega_{2}})
\]
converges in $\mu_1$-probability on $\Omega_{1}$ as $n\to\infty$,
and the function $\phi:\Omega_{1}\to\mathbb{R\cup\{-\infty\}}$ is measurable with $\phi>-\infty$ on a set of positive $\mu_{1}$-probability.
Here we possibly have to pass to a subsequence to guarantee that there is $ \widetilde{\omega_2} $ where $ \xi_n(.,\widetilde{\omega_2}) $ 
converges in probability. 
We assume from now on that we have already passed to such a subsequence.
We then define $r_{n}:=\psi_{n}(\widetilde{\omega_{2}})$, for
$n\geq1$. Furthermore we obtain that the sequence
\[
\psi_{n}-r_{n}=\xi_{n}-\phi_{n}-r_{n}
\]
converges in $\mu_2$-probability on $\Omega_{2}$ as $n\to\infty$
to a function $\psi:\Omega_{2}\to\mathbb{R\cup\{-\infty\}}$. Indeed, the right hand side does not depend on the first variable $\omega_{1}$ by construction. We
choose a measurable subset $A_{1}\subset\Omega_{1}$ of positive $\mu_{1}
$-probability with $\phi|_{A_{1}}>-\infty$. The sequence $(\xi_{n})_{n\geq1}$
converges in probability to $\xi$ and the sequence $(\phi_{n}+r_{n})_{n\geq1}$
converges in probability to $\phi$ on $A_{1}$, so the difference converges
in probability on $A_{1}\times\Omega_{2}$ to a limit $\psi$, which does not
depend on the first variable. Hence in particular
\[
\xi(\omega_{1},\omega_{2})=\phi(\omega_{1})+\psi(\omega_{2})
\]
$\mu_{1}\otimes\mu_{2}$-almost sure.

Now we assume that $\xi$ is almost surely equal to $-\infty$. We know that,
for each fixed $0<\alpha<1$, the $\alpha$-quantiles
\begin{align*}
q_{n}^{1}(\alpha) &  :=\inf\{x;~\mu_{1}(\phi_{n}\geq x)\leq\alpha\},\\
q_{n}^{2}(\alpha) &  :=\inf\{y;~\mu_{2}(\psi_{n}\geq y)\leq\alpha\}
\end{align*}
have the property
\[
\lim_{n\to\infty}(q_{n}^{1}(\alpha)+q_{n}^{2}(\alpha))=-\infty
\]
by assumption. Indeed, otherwise there is a subsequence $(n_{k})_{k\geq1}$ and
a lower bound $M>-\infty$ such that $q_{n_{k}}^{1}(\alpha)+q_{n_{k}}
^{2}(\alpha)\geq M$ for all $k\geq1$. Hence $\phi_{n_{k}}+\psi_{n_{k}}\geq M
$ with $\mu_{1}\otimes\mu_{2}$-probability greater than $\alpha^{2}$
for $k\geq1$. However, this contradicts the convergence in probability to
$-\infty$. By choosing a sequence $(\alpha_{n})_{n\geq1}$ with
$\alpha_{n}\downarrow0$
slowly enough as $n\to\infty$, we can maintain this property, i.e.
\[
\lim_{n\to\infty}(q_{n}^{1}(\alpha_{n})+q_{n}^{2}(\alpha_{n}
))=-\infty.
\]
We define now
\begin{align*}
z_{n} &  :=\frac{1}{2}(q_{n}^{1}(\alpha_{n})+q_{n}^{2}(\alpha_{n}))\\
r_{n} &  :=\frac{1}{2}(q_{n}^{2}(\alpha_{n})-q_{n}^{1}(\alpha_{n}))\\
&  =q_{n}^{2}(\alpha_{n})-z_{n}=-(q_{n}^{1}(\alpha_{n})-z_{n}).
\end{align*}
We have $\lim_{n\to\infty}z_{n}=-\infty$ and hence $\phi_{n}
+r_{n}\to-\infty$ and $\psi_{n}-r_{n}\to-\infty$ converges in the respective probabilities as $n\to\infty$ by construction.
\end{proof}

In the proof of the subsequent proposition we shall apply a Komlos-type
result, which we cite here from \cite[Lemma A1.1]{delsch:94} for the sake of completeness.

\begin{lemma}
\label{komlos} Let $(f_{n})_{n\geq1}$ be a sequence of $\mathbb{R}_{\geq0}%
\cup\{+\infty\}$-valued random variables on a probability space $(\Omega
,\mathcal{F},P)$. Then there is a sequence $g_{n}\in\operatorname*{conv}%
(f_{n},f_{n+1},\dots)$, the convex hull of $f_{n},f_{n+1},\dots$, for $n\geq
1$, which converges almost surely to a random variable $g$ taking values in
$\mathbb{R}_{\geq0}\cup\{+\infty\}$.
\end{lemma}

The interesting feature of this result is that we do not need any
integrability assumption on the sequence $(f_{n})_{n\geq1}$ to obtain an
almost sure convergence of convex combinations. It suffices to assume
non-negativity of the functions $f_{n}$. This situation is related to the
setting of \cite{brasch:99}, where it is shown that a Bipolar-Theorem can be
formulated for the non-locally-convex space $L^{0} $, provided that one
restricts to positive functions.

\begin{proposition}
Let $c:X\times Y\to\mathbb{R}_{\geq0}\cup\{+\infty\}$ be a lower
semi-continuous cost function on polish spaces $X,Y$. Given Borel probability
measures $\mu$ on $X$, $\nu$ on $Y$ and an optimizer $\widehat{\pi}\in\Pi
(\mu,\nu)$ with $I(\widehat{\pi})<\infty$. We assume that $c$ is $\mu
\otimes\nu$-almost surely finitely valued. Then there exist Borel functions
$\phi,\psi$ on $X,Y$, respectively, taking values in $\mathbb{R}\cup
\{-\infty\}$, such that
\[
\phi(x)+\psi(y)\leq c(x,y)
\]
for $x\in X$ and $y\in Y$, and $\phi(x)+\psi(y)=c(x,y)$ almost surely with
respect to $\widehat{\pi}$. In particular the transport plan $\widehat{\pi}$
is strongly $c$-monotone.\label{measure theoretic}
\end{proposition}

\begin{proof}
By \cite{ambpra:02}, Th. 3.1, we know that there is a maximizing sequence
$(\phi_{n},\psi_{n})$ of bounded Borel functions on $X,Y$, respectively, such
that
\[
\lim_{n\to\infty}E_{\widehat{\pi}}(\xi_{n})=\lim_{n\to\infty
}E_{\mu\otimes\nu}(\xi_{n})=E_{\widehat{\pi}}(c)=I(\widehat{\pi}),
\]
where $\xi_{n}(x,y)=\phi_{n}(x)+\psi_{n}(y)\leq c(x,y)$ for all $(x,y)\in
X\times Y$ and the convergence of $(E_{\mu\otimes\nu}(\xi_{n}))_{n\geq0}$ is
monotonely increasing. By passing to convex combinations (apply the above
Lemma \ref{komlos} with $f_{n}=c-\xi_{n}$ for $n\geq0$) we may assume that
$\xi_{n}$ converges $\mu\otimes\nu$-almost surely to a Borel function $\xi$
taking values in $\mathbb{R}\cup\{-\infty\}$. By Lemma \ref{convergence} we
know that
\[
\xi(x,y)=\phi(x)+\psi(y),
\]
with functions $\phi$ and $\psi$, which are $\mu$- and $\nu$-almost surely
defined, respectively. Furthermore -- by passing to a subsequence if necessary
-- there is a sequence of real numbers $(r_{n})_{n\geq1}$ such that $\phi
_{n}+r_{n}\to\phi$ and $\psi_{n}-r_{n}\to\psi$ as
$n\to\infty$ almost surely with respect to $\mu$ and $\nu$,
respectively. Choosing appropriate nullsets in $X$ and $Y$, we can redefine
$\phi,\phi_{n}$ and $\psi$,$\psi_{n}$ on these nullsets by $-\infty$, such
that $\phi_{n}+r_{n}\to\phi$ and $\psi_{n}-r_{n}\to\psi$
\emph{surely}, without violating the now \emph{sure inequality}
\[
\xi(x,y)=\phi(x)+\psi(y)\leq c(x,y).
\]
Obviously $E_{\widehat{\pi}}(|c-\xi_{n}|)=E_{\widehat{\pi}}(c-\xi
_{n})\to0$, as $n\to\infty$, by assumption, hence
$\phi(x)+\psi(y)=c(x,y)$ almost surely with respect to $\widehat{\pi}$.
\end{proof}

Our main Theorem below states that strongly $c$-monotone transport plans are
always optimal. For sake of clarity we shall formulate an elementary lemma on
monotone convergence of truncations, whose proof is obvious.

\begin{lemma}
Let $a,b$ be real numbers. We define
\begin{align*}
a_{n}  &  :=(-n\lor a)\land n,\\
b_{n}  &  :=(-n\lor b)\land n,\\
\xi_{n}  &  :=a_{n}+b_{n}.
\end{align*}
Then we have $\xi_{0}=0$ and $(\xi_{n})_{n\geq0}$ converges monotonically to
$a+b$, i.e. if $a+b\geq0$ the sequence $(\xi_{n})_{n\geq0}$ increases to $a+b$
and if $a+b\leq0$ the sequence decreases to $a+b$.
\end{lemma}

\begin{theorem}
Let $c:X\times Y\to\mathbb{R}_{\geq0}\cup\{+\infty\}$ be a lower
semi-continuous cost function on polish spaces $X,Y$. Given Borel probability
measures $\mu$ on $X$, $\nu$ on $Y$ and a strongly $c$-monotone transport plan
$\pi\in\Pi(\mu,\nu)$, then $\pi$ is an optimizer of the Monge-Kantorovich problem.\label{characterization}
\end{theorem}

\begin{proof}
We assume that there is a transport plan $\pi_{0}$ with finite cost $I(\pi
_{0})<\infty$. We aim to show that $I(\pi)\leq I(\pi_{0})$. By the assumption
of strong $c$-monotonicity of $\pi$ there exist Borel-measurable functions
$\phi,\psi$ on $X,Y$, respectively, taking values in $\mathbb{R}\cup
\{-\infty\}$, such that for all $(x,y)\in X\times Y$
\[
\phi(x)+\psi(y)\leq c(x,y),
\]
where equality holds $\pi$-almost surely. We define
\begin{align*}
\phi_{n}(x)  &  :=(-n\lor\phi(x))\land n,\\
\psi_{n}(y)  &  :=(-n\lor\psi(y))\land n,\\
\xi_{n}(x,y)  &  =\phi_{n}(x)+\psi_{n}(y),\\
\xi(x,y)  &  =\phi(x)+\psi(y),
\end{align*}
for $(x,y)\in X\times Y$ and $n\geq0$. Due to the previous Lemma observe that
$\xi_{n}\uparrow\xi$ on $\{\xi\geq0\}$ and $\xi_{n}\downarrow\xi$ on
$\{\xi\leq0\}$, as $n\to\infty$.

Additionally, $E_{\pi}(\xi)=E_{\pi}(c)$ exists, taking possibly the value
$+\infty$, since equality $\xi=c$ holds $\pi$-almost surely and $c\geq0$. The
integral of $\xi$ with respect to $\pi_{0}$ exists, too, as
\[
E_{\pi_{0}}(\xi)\leq E_{\pi_{0}}(c)<\infty.
\]
Note that $E_{\pi_{0}}(\xi)$ possibly equals $-\infty$. By the assumption on
equal marginals of $\pi$ and $\pi_{0}$ we obtain
\begin{align*}
E_{\pi}(\xi_{n})  &  =E_{\pi}(\phi_{n})+E_{\pi}(\psi_{n})\\
&  =E_{\pi_{0}}(\phi_{n})+E_{\pi_{0}}(\psi_{n})\\
&  =E_{\pi_{0}}(\xi_{n}),
\end{align*}
for $n\geq0$, hence
\[
E_{\pi}(\xi_{n}1_{\{\xi\geq0\}}+\xi_{n}1_{\{\xi\leq0\}})=E_{\pi_{0}}(\xi
_{n}1_{\{\xi\geq0\}}+\xi_{n}1_{\{\xi\leq0\}}).
\]
By our previous considerations we can pass to the limits and obtain $E_{\pi
}(\xi)=E_{\pi_{0}}(\xi)$. Indeed the limits are monotone on $\{\xi\geq0\}$ and
$\{\xi\leq0\}$, the expectation with respect to $\pi$ of $\xi_{-}$ is finite,
namely $E_{\pi}(\xi_{-})=0$, and the expectation with respect to $\pi_{0}$ of
$\xi_{+}$ is finite, $E_{\pi_{0}}(\xi_{+})<\infty$. Hence the limits of
$E_{\pi}(\xi_{n})=E_{\pi_{0}}(\xi_{n})$ exist as $n\to\infty$ and are
equal. Consequently $I(\pi)=E_{\pi}(\xi)=E_{\pi_{0}}(\xi)\leq I(\pi_{0})$.
\end{proof}

The following Corollary answers Problem 2.25 of \cite{vil:03} pertaining to
the quadratic cost function affirmatively. In fact, all finite, lower
semi-continuous cost functions are covered and no additional assumptions on
integrability or measurability are necessary.

\begin{corollary}
Let $c:X\times Y\to\mathbb{R}_{\geq0}$ be a finitely valued, lower
semi-continuous cost function on polish spaces $X,Y$. Given Borel probability
measures $\mu$ on $X$, $\nu$ on $Y$ and a $c$-monotone transport plan $\pi
\in\Pi(\mu,\nu)$, then $\pi$ is an optimizer of the Monge-Kantorovich problem.
\end{corollary}

\begin{proof}
Due to Proposition \ref{existence}, for a finitely valued, lower
semi-continuous cost function $c$, a $c$-monotone transport plan $\pi$ is
strongly $c$-monotone. Hence Theorem \ref{characterization} applies.
\end{proof}

The subsequent example, which is taken from \cite{ambpra:02}, Example 3.5,
shows that there are $c$-monotone transport plans, which are not strongly $c$-monotone.

\begin{example}
This example is just a re-interpretation of the illuminating Example
3.5. in \cite{ambpra:02}. Take $X=Y=T^{1}=\{\exp(2\pi ir), \, 0 \leq r
< 1 \}\subset\mathbb{C}$ the $1$-dimensional Torus with uniform
distribution $\mu=\nu$ and choose $\alpha\in T^{1}$, so that
$\{\alpha^{n},n\geq0\}$ is dense in $T^{1}$, i.e. $\alpha=\exp(2\pi
ir)$ with $r$ an irrational number: we define a lower semicontinuous
cost function $c:X\times Y\to
\mathbb{R}_{\geq0}\cup\{+\infty\}$ via
\[
c(x,y)=\left\{
\begin{array}
[c]{c}%
1\text{ if }x=y\\
2\text{ if }x=y\alpha
\end{array}
\right.
\]
and $+\infty$ elsewhere. Then there are two disjoint $c$-monotone subsets of
$X\times Y$, namely
\begin{align*}
\Gamma_{1}  &  =\{(x,x),x\in T^{1}\},\\
\Gamma_{2}  &  =\{(x,x\alpha),x\in T^{1}\}.
\end{align*}
For the first set strong $c$-monotonicity is clear with $\phi_{1}(\exp(2\pi
ir))=r$ for $r\in\lbrack0,1[$, and $\psi_{1}(\exp(ir))=\inf_{x}(c(x,\exp
(ir))-\phi_{1}(x))=1-r$ for $r\in\lbrack0,1[$. Then we obtain the result, that
$\phi_{1}(x)+\psi_{1}(y)\leq c(x,y)$ with equality if $x=y$.

The second set $\Gamma_{2}$ is $c$-monotone (see \cite{ambpra:02}, Example
3.5.), but not strongly $c$-monotone (this follows from non-optimality, but we
shall also show it directly).

Take $n\geq1$ and $(x_{i},x_{i}\alpha)$ for $i=1,\dots,n$. Then for all
$\sigma\in\mathfrak{S}_{n}$,
\[
\sum_{i=1}^{n}c(x_{i},x_{i}\alpha)\leq\sum_{i=1}^{n}c(x_{i},x_{\sigma
(i)}\alpha).
\]
Indeed, otherwise there is $\sigma\in\mathfrak{S}_{n}$ and at least one
$x_{i}\in T^{1}$ such that $c(x_{i},x_{\sigma(i)}\alpha)=1$ and therefore
\[
x_{i}=x_{\sigma(i)}\alpha
\]
Iteration of this equation yields $\alpha^{m}=1$ for some $m\geq1$, a
contradiction showing the $c$-monotonicity of $\Gamma_{2}$.

We shall now show that $\Gamma_{2}$ fails to be strongly $c$-monotone. In
fact, we shall show the stronger assertion that, for every Borel measure
$B\subset T^{1}$ with $\mu(B)>0$ the set
\[
\Gamma_{3}=\{(x,x\alpha),x\in B\}
\]
fails to be strongly $c$-monotone. Indeed, we suppose that there are
Borel-measurable $\phi_{2}$ and $\psi_{2}$ taking value in $\mathbb{R}%
\cup\{-\infty\}$, such that
\begin{equation}
\phi_{2}(x)+\psi_{2}(y)\leq c(x,y) \tag{2}\label{e1}
\end{equation}
for all $(x,y)\in T^{1}\times T^{1}$, with equality holding for $x\alpha=y$
and $x\in B$. Using the special form of $c$, we see that (\ref{e1}) implies
\begin{align}
\phi_{2}(z)+\psi_{2}(z)  &  \leq1,\tag{3}\label{e2}\\
\phi_{2}(x)+\psi_{2}(x\alpha)  &  =2, \tag{4}\label{e3}
\end{align}
for $x\in B$ and $z\in T^{1}$. Combining (\ref{e2}) and (\ref{e3}) we obtain
\begin{equation}
\phi_{2}(x)\geq\phi_{2}(x\alpha)+1 \tag{5}\label{e4}
\end{equation}
for $x\in B$ or for $x\alpha\in B$. Defining $B_{k}:=B\cap\phi_{2}
^{-1}([k,k+1[)$ for $k\in\mathbb{Z}$, we obtain a partition $(B_{k}
)_{k\in\mathbb{Z}}$ of $B$. We consider now an equivalence relation on $T^{1}
$, namely $x\sim_{\alpha}y$ if $x\alpha^{m}=y$ for some $m\in\mathbb{Z}$. Fix
$k\in\mathbb{Z}$; due to \ref{e4} there can be at most one element of each
equivalence class $[x]\subset T^{1}$ in $B\cap\phi_{2}^{-1}([k,k+1[)$. On the
other hand every Borel set $A\subset T^{1}$, which contains at most one
element of each equivalence class $[x]$, has measure $\mu(A)=0$ (since then
$\cup_{m\in\mathbb{Z}}(A\alpha^{m})\subset T^{1}$, $A\alpha^{m}$ are pairwise
disjoint and $\mu(A\alpha^{m})=\mu(A)$ for $m\in\mathbb{Z}$ by translation
invariance). Consequently we obtain a contradiction to the assumption that the
sets $B_{k}$ are Borel sets.
\end{example}

We can finally formulate an approximation result for the optimal transport
plan of Monge-Kantorovich problems. The result can also be interpreted as
discretization result for numerical purposes. In particular, we obtain --
without any moment conditions -- that the adherence points of candidate
approximations are optimizers, provided the cost function $c$ is continuous.

\begin{definition}
Let $c:X\times Y\to\mathbb{R}_{\geq0}\cup\{+ \infty\}$ be a
lower semi-continuous cost function on polish spaces $X,Y$. Given Borel
probability measures $\mu$ on $X$ and $\nu$ on $Y$. Assume that there
are sequences $(\mu_{n})_{n\geq1}$ and $(\nu_{n})_{n\geq1}$ of Borel
probability measures converging weakly to $\mu$ and $\nu$ on the
respective spaces $X$ and $Y$. We denote by $\pi_{n}$ a solution of the
Monge-Kantorovich problem associated to
$(\mu_{n},\nu_{n},c)$ for $n\geq1$ and we assume that
\[
\int_{X\times Y}c(x,y)\pi_{n}(dx,dy)<\infty,
\]
for each $n\geq1$. We call the sequence $(\pi_{n})_{n\geq1}$ an
\emph{approximating sequence of optimizers}.
\end{definition}

\begin{theorem}
Let $c:X\times Y\to\mathbb{R}_{\geq0}$ be a finitely valued,
continuous cost function on polish spaces $X,Y$. Let
$(\pi_{n})_{n\geq1}$ be an approximating sequence of optimizers
associated to weakly converging sequences $\mu_{n}\to\mu$ and
$\nu_{n}\to\nu$ as $n\to\infty$. Then there is a
subsequence $(\pi_{n_{k}})_{k\geq0}$ converging weakly to a transport
plan $\pi$ on $X\times Y$, which optimizes the Monge-Kantorovich
problem $(\mu,\nu,c)$. Any other converging subsequence of
$(\pi_{n})_{n\geq0}$ also converges to an optimizer of the
Monge-Kantorovich problem, i.e. the non-empty set of adherence points
of $(\pi_{n})_{n\geq0}$ is a set of optimizers.\label{discretization}
\end{theorem}

\begin{proof}
Fix $\epsilon>0$. By Prohorov's Theorem (see for instance \cite{kal:01}, Th.
16.3) there are compact sets $K\subset X$ and $L\subset Y$ such that $\mu
_{n}(X\setminus K)\leq\epsilon$ and $\nu_{n}(Y\setminus L)\leq\epsilon$ for
all $n$. Hence $\pi_{n}(X\times Y\setminus K\times L)\leq\pi_{n}
(X\times(Y\setminus L))+\pi_{n}((X\setminus K)\times Y)\leq2\epsilon$. Again
by Prohorov's Theorem we know that there is a weakly converging subsequence
$(\pi_{n_{k}})_{k\geq0}$ with weak limit $\pi$. We have to prove that $\pi$ is
$c$-monotone, since by Proposition \ref{existence} the transport plan $\pi$ is
then strongly $c$-monotone and hence optimal. We know that for $m\geq1$ and
points $(x_{i},y_{i})\in\operatorname*{supp}\pi$ for $i=1,\dots,m$, there
exist sequences of points $(x_{i}^{k},y_{i}^{k})_{k\geq0}$ such that
$(x_{i}^{k},y_{i}^{k})\in\operatorname*{supp}\pi_{n_{k}}$ and $(x_{i}
^{k},y_{i}^{k})\to(x_{i},y_{i})$ as $k\to\infty$ (see again
\cite{kal:01}, Th. 16.3). Since $\pi_{n_{k}}$ is $c$-monotone by Proposition
\ref{measure theoretic} and Remark \ref{strong implies usual}, and since $c$
is continuous, $\operatorname*{supp}\pi_{n_{k}}$ is a $c$-monotone subset of
$X\times Y$, hence
\[
\sum_{i=1}^{m}c(x_{i}^{k},y_{i}^{k})\leq\sum_{i=1}^{m}c(x_{\sigma(i)}
^{k},y_{\sigma(i)}^{k})
\]
for permutations $\sigma\in\mathfrak{S}_{m}$, whence also for
$k\to \infty$ by continuity of $c$. By the same reasoning any
converging subsequence of $(\pi_{n})_{n\geq0}$ converges to a
$c$-monotone transport plan, which is by Proposition \ref{existence} a
strongly $c$-monotone one and hence optimal by Theorem
\ref{characterization}.
\end{proof}

\begin{remark}
In particular the set of optimizers of a Monge-Kantorovich problem is compact in the weak topology on 
probability measures, since the constant sequences $ \mu_n = \mu $ and $ \nu_n = \nu $ for $ n \geq 0 $ are approximating sequences. 
This also yields that we cannot improve the result to the assertion on convergence of the sequence, but we point out that the convergence of the sequence in Theorem \ref{discretization} holds true if and only if the optimizer is unique.
\end{remark}

\begin{example}
Let $c$ be a continuous, finitely valued cost function on the product $X\times
Y$ of two polish spaces $X,Y$ and let $\mu,\nu$ be probability measures
thereon. By the law of large numbers we can always find empirical
distributions (see for instance \cite{kal:01}, Prop. 4.24)
\begin{align*}
\mu_{n} &  =\frac{1}{n}\sum_{i=1}^{n}\delta_{x_{i}},\\
\nu_{n} &  =\frac{1}{n}\sum_{i=1}^{n}\delta_{y_{i}},
\end{align*}
which converge weakly to $\mu$ and $\nu$, respectively. Each associated
approximating sequence of optimizers $(\pi_{n})_{n\geq1}$ is given
through a sequence of permutations $\sigma_{n}\in\mathfrak{S}_{n}$,
$n\geq1$,
and%
\[
\pi_{n}=\frac{1}{n}\sum_{i=1}^{n}\delta_{(x_{i},y_{\sigma_{n}(i)})}.
\]
Notice that always $I(\pi_{n})=\int_{X\times
Y}c(x,y)\pi_{n}(dx,dy)<\infty$ for each fixed $n\in\mathbb{N}$. We
deduce from Theorem \ref{discretization} that all adherence points of
$(\pi_{n})_{n\geq1}$ are optimizers of the Monge-Kantorovich problem
and that there is at least one adherence point. Notice that each $
\pi_n $ for $ n \geq 1 $ is obtained by solving a finite optimization
problem.
\end{example}

The authors would like to thank the participants of the Seminar
''Topics in Optimal Transportation'' for the fruitful discussions on C.
Villani's book. The Seminar is supported by the FWF-Wissenschaftskolleg
W 800.

It is a particular pleasure to thank Lugi Ambrosio, Aldo Pratelli and
Cedric Villani for their very helpful comments on a previous version of
this note. After completion of our work we were kindly informed about
the interesting and independent work \cite{pra:06} on the same subject
in a slightly different framework. We also thank the anonymous referee
for her/his pertinent and insightful report, which led to an improvement of
the paper.

\end{document}